\documentclass[a4paper,12pt, draft]{article}
\usepackage{a4wide}
\usepackage{amsmath,amsthm,amssymb}
\usepackage[mathscr]{eucal}
\theoremstyle{plain}
\newtheorem{proposition}{Proposition}
\newtheorem{theorem}{Theorem}
\newtheorem{lemma}{Lemma}
\theoremstyle{definition}
\newtheorem{df}{Definition}
\theoremstyle{remark}
\newtheorem{corollary}{Corollary}
\newcommand{\dom}{\mathop{\rm dom}\nolimits}
\newcommand{\ran}{\mathop{\rm ran}\nolimits}
\newcommand{\rank}{\mathop{\rm rank}\nolimits}

\begin{document}
\author{Galyna Tsyaputa}
\title{Deformed multiplication in the semigroup $\mathcal{PT}_n$}
\date{}
\maketitle
\begin{abstract}
Pairwise non isomorphic semigroups obtained from the semigroup
$\mathcal{PT}_n$ of all partial transformations by the deformed
multiplication proposed by Ljapin are classified.
\end{abstract}

\section{Introduction}
Let $X$ and $Y$ be two nonempty sets, $S$ be the set of maps from
$X$ to $Y$. Fix some $\alpha : Y\rightarrow X$ and define the
multiplication of elements in $S$ in the following way: $\phi\circ
\psi = \phi\alpha\psi$ (the composition of the maps is from left
to right). The action defined by this rule is associative. In his
famous monograph ~\cite[p.353]{Ljap} Ljapin set the problem of
investigation of the properties of this semigroup depending on the
restrictions on set $S$ and map $\alpha$.

 Magill \cite{Mag} studies this problem in the case when $X$ and $Y$ are
 topological spaces and the maps are continuous.
 In particular under the assumption that $\alpha$ be onto he describes the automorphisms
 of such semigroups and determines the isomorphism criterion.

A bit later Sullivan \cite{Sul} proves that if $|Y|\leq|X|$ then
Ljapin's semigroup is embedding into transformation semigroup on
the set $X\cup \{a\}$, $a\notin X$.

For us the important case is if $X=Y$, $T_X$ is a transformation
semigroup on the set $X$, $\alpha\in T_X$. Symons \cite{Sym}
establishes the isomorphism criterion for such semigroups and
investigates the properties of their automorphisms.

The latter problem may be generalized to arbitrary semigroup $S:$
for a fixed $a\in S$ define the operation $*_a$ via $x*_a y=xay$.
 The obtained semigroup is denoted by $(S,*_a)$ and operation $*_a$ is
called the multiplication deformed by element $a$ (or just the
deformed multiplication).

In the paper we classify with respect to isomorphism all
semigroups which are obtained from the semigroup of all partial
transformations of an $n-$element set by the deformed
multiplication.

\section{$\mathcal{PT}_n$ with the deformed multiplication}

Let $\mathcal{PT}_n $ be the semigroup of all partial
transformations of the set $N=\{1, 2,\ldots , n\}$. For arbitrary
partial transformation $a\in\mathcal{PT}_n $ denote by $\dom (a)$
the domain of $a$, and by $\ran (a)$ its image. The value $|\ran
(a)|$ is called the rank of $a$ and is denoted by $\rank (a)$.
Denote $Z(a)=N\setminus \dom (a)$ the set of those elements on
which transformation $a$ is not defined, and denote $z_a $ the
cardinality of this set. The type of $a\in\mathcal{PT}_n $ we call
the set $(\alpha_1 ,\alpha_2 ,\dots ,\alpha_n)$, where $\alpha_k $
is the number of those elements $y\in N$, the full inverse image
$a^{-1} (y)$ of which contains exactly $k$ elements. It is obvious
that $1\cdot\alpha_1 +2\cdot\alpha_2 +\dots +n\cdot\alpha_n
=|\dom(a)|$, and the sum $\alpha_1 +\alpha_2 +\dots +\alpha_n$ is
equal to the rank of $a$.

\begin{df} Element $x\in\mathcal{S}$ is called left (right) annihilator
of semigroup $\mathcal{S}^0$ provided that $xs=0$ ($sx=0$),
$s\in\mathcal{S}$.

Element which is both the left and the right annihilator is called
the annihilator of semigroup $\mathcal{S}$.
\end{df}
\begin{proposition}
\label{tvr1} Transformations $x\in\mathcal{PT}_n$ such that $\ran
(x)\subset Z(a)$ are left annihilators of semigroup
$(\mathcal{PT}_n ,*_a)$. The number of left annihilators equals
$(z_a +1)^n$.

Transformations $y\in\mathcal{PT}_n$ such that $Z(y)\supset
\ran(a)$ are right annihilators of semigroup $(\mathcal{PT}_n
,*_a)$. The number of right annihilators equals $(n+1)^{n-\rank
(a)}$.

Transformations $c\in\mathcal{PT}_n$ such that $\ran (c)\subset
Z(a)$ and $Z(c)\supset\ran (a)$ are annihilators of
$(\mathcal{PT}_n ,*_a)$, moreover, the number of annihilators
equals $(z_a +1)^{n-\rank (a)}$.

\end{proposition}
\begin{proof}Let $x\in\mathcal{PT}_n $ be left annihilator of $(\mathcal{PT}_n
,*_a)$. This means that for arbitrary $u$ from $\mathcal{PT}_n $
$x*_a u=0$ or, what is the same, $xa\cdot u=0 $. Therefore $xa$ is
left zero of semigroup $\mathcal{PT}_n $, that is, a nowhere
defined map. The latter is possible if and only if $\ran
(x)\subset Z(a)$. Evidently, the condition $\ran (x)\subset Z(a)$
is sufficient for partial transformation $x$ be the left
annihilator of $(\mathcal{PT}_n ,*_a)$. Now it is clear, that the
number of such transformations is equal to $(z_a +1)^n $.

Let $y\in\mathcal{PT}_n $ be the right annihilator of semigroup
$(\mathcal{PT}_n ,*_a)$. Then for any $v$ from $\mathcal{PT}_n $
$v\cdot ay=0$ and $ay$ is right zero of semigroup $\mathcal{PT}_n
$, that is, nowhere defined map. This is equivalent to
$Z(y)\supset\ran (a)$. Therefore right annihilators of semigroup
$(\mathcal{PT}_n ,*_a) $ are those partial transformations
$\mathcal{PT}_n $ which are defined only in elements
$N\setminus\ran (a)$. It is clear that the number of such
transformations equals $(n+1)^{n-\rank (a)}$.

Now the statement about annihilators follows from the definition
and the above arguments.
\end{proof}
On semigroup $(\mathcal{PT}_n ,*_a)$  define the equivalence
relation $\sim _a $ by the rule: $x\sim _a y$ if and only if $x*_a
u=y*_a u$ for all $u\in \mathcal{PT}_n$. Analogously on
$(\mathcal{PT}_n ,*_b)$ define the relation $\sim _b$.
\begin{lemma}
\label{lemma1} For an arbitrary isomorphism $\varphi :
(\mathcal{PT}_n ,*_a)\rightarrow (\mathcal{PT}_n ,*_b)$ $\varphi
(x)\sim _b\varphi (y)$ if and only if $x\sim _a y$.
\end{lemma}
\begin{proof} In fact, let $\varphi : (\mathcal{PT}_n
,*_a)\rightarrow (\mathcal{PT}_n ,*_b)$ be isomorphism and $x\sim
_a y$. Then for all $u\in\mathcal{PT}_n$: $x*_a u=y*_a u$,
therefore $\varphi (x)*_b\varphi(u)=\varphi (y)*_b\varphi(u)$.
However $\varphi (u)$ runs over the whole set $\mathcal{PT}_n$,
hence $\varphi (x)\sim_b\varphi (y)$. Since the inverse map
$\varphi^{-1} : (\mathcal{PT}_n ,*_b)\rightarrow (\mathcal{PT}_n
,*_a)$ is also isomorphism, $\varphi (x)\sim_b\varphi (y)$ implies
$x\sim _a y$. Therefore, $x\sim_a y$ if and only if $\varphi
(x)\sim_b \varphi (y)$.
\end{proof}
\begin{lemma}
$x\sim_a y$ if and only if $xa=ya$.
\end{lemma}
\begin{proof}
Obviously, the equality $xa=ya$ implies $x\sim_a y$. Now let
$xa\neq ya$. Then there exists $k$ in $N$, such that $(xa)(k)\neq
(ya)(k)$. Chose element $u$ in $\mathcal{PT}_n $ which has
different images in the points $(xa)(k)$ and $(ya)(k)$. Then $x*_a
u=xau$ and $y*_a u=yau$ have different images in $k$. Hence $x*_a
u\neq y*_a u$ and $x\nsim _a y$.
\end{proof}
\begin{theorem}
Semigroups $(\mathcal{PT}_n ,*_a)$ and $(\mathcal{PT}_n ,*_b)$ are
isomorphic if and only if partial transformations $a$ and $b$ have
the same type.
\end{theorem}
\begin{proof}\textsl{Necessity.} Let $(\mathcal{PT}_n ,*_a)$ and $(\mathcal{PT}_n ,*_b)$
 be isomorphic. By lemma~\ref{lemma1} arbitrary isomorphism between $(\mathcal{PT}_n ,*_a)$
 and $(\mathcal{PT}_n ,*_b)$ maps equivalence classes of the relation $\sim_a$
 into equivalence classes of the relation $\sim_b $.
 Therefore for equivalence relations $\sim_a$ and $\sim_b$ cardinalities and the number of equivalence
 classes must be equal. We show that by the cardinalities of the classes of the relation $\sim_a$
 the type $(\alpha_1,\alpha_2,\dots ,\alpha_n)$ of transformation $a$ can
be found uniquely.

Denote $\rho_a$ the partition of set $\{1,2,\dots ,n\}$ induced by
$a\in\mathcal{PT}_n$ (that is, $l$ and $m$ belong to the same
block of the partition $\rho_a$ provided that $a(l)=a(m)$; $Z(a)$
forms a separate block). We count the cardinality of equivalence
class $\overline{x_0}=\{x|xa=x_0 a\}$ of the relation $\sim_a$ for
the fixed transformation $x_0\in\mathcal{PT}_n$. Consider element
\begin{displaymath}y:=x_0 a=\left(
\begin{array}{lccccccccr}
 i_1 &i_2& \dots & i_p & i_{p+1} & \dots & i_n \\
 y_1 & y_2 & \dots & y_p & \varnothing &\dots &\varnothing \end{array}
 \right).
\end{displaymath}
 Obviously $y_i$ belongs to the image of $a$, $i=1,\dots , p$.
 Denote $N_a (a_i)$ the block of the partition $\rho_a $ which is defined by $a_i$ in the image of
 $a$. Denote $n_a (a_i) $ the cardinality of this block. The equality $xa=y$
 is equivalent to that for every $i$
$(xa)(i)=y_i$, or, what is the same, $x(i)\in N_a (y_i)$. Hence
$x(i)$ can be chosen in $n_a (y_i)$ ways, $i=1,\dots ,p$. For
$i=p+1,\dots ,n$ the meaning of $x(i)$ can be chosen in $z_a +1$
ways. Since the images of $x$ in different points are chosen
independently, partial transformation $x$ can be defined in
\begin{equation}
\label{odynPTn} n_a (y_1)n_a(y_2)\dots n_a(y_p)(z_a +1)^{n-p}
\end{equation}
ways, and this gives the cardinality of class $\overline{x_0}$.

It is clear that $\overline{0}$ is the set of all left
annihilators of the semigroup $(\mathcal{PT}_n ,*_a)$. By
proposition~\ref{tvr1} the cardinality of class $\overline{0}$
equals $(z_a +1)^n$. Hence the value $z_a$ is defined by abstract
property the semigroup $(\mathcal{PT}_n ,*_a)$. Denote $m$ the
smallest cardinality of the blocks of the partition $\rho_a$. Then
the cardinality of the equivalence class of the relation $\sim_a$
is the least if in~(\ref{odynPTn}) all multipliers are equal to
$\min (m, z_a +1)$. Consider the corresponding cases.

Let $m>z_a +1$. Then class $\overline{0}$ is the only equivalence
class of the least possible cardinality. The next larger class
contains $m(z_a +1)^{n-1}$ transformations. Now count the number
of different equivalence classes of the relation $\sim_a$ which
have the cardinality $m(z_a +1)^{n-1}$. Since $m>z_a +1$ and this
is the least of the cardinalities of the blocks of the partition
$\rho_a$, to make the product in~(\ref{odynPTn}) be equal to
$m(z_a +1)^{n-1}$ there should be $p=1$ and $n_a (y_1)=m$. Element
$i_1$ can be chosen in $n$ ways. We know that $|\{t:n_a
(t)=m\}|=\alpha_m$, so $y_1$ can have $\alpha_m$ different
meanings. Therefore, there are $n\cdot\alpha_m $ different
transformations $y=x_0 a$, such that the corresponding class
$\overline{x_0}$ contains $m(z_a +1)^{n-1}$ elements. Therefore by
the relation $\sim_a$ we may define the number $m$ and the meaning
$\alpha_m$ of the first nonzero component of the type $(\alpha_1
,\alpha_2,\dots,\alpha_n)$ of element $a$. Components $\alpha_l$
for $l>m$ can be defined recursively. Assume that $\alpha_1
,\alpha_2,\dots,\alpha_l$ are found. For the relation $\sim_a$
denote $C$ the number of equivalence classes of the cardinality
$l\cdot (z_a +1)^{n-1}$. Then $C$ is equal to the number of sets
$(p; i_1, i_2,\dots ,i_p),$ where some of $i_1, i_2,\dots ,i_p$
may coincide in general, such that
\begin{equation}
\label{dvaPTn} n_a (y_{i_1})n_a(y_{i_2})\dots n_a(y_{i_p})(z_a
+1)^{n-p}=l\cdot(z_a +1)^{n-1}.
\end{equation}
Since $\alpha_1,\alpha_2,\dots ,\alpha_{l-1}$ are known, we may
find the number $A$ of sets $(p; i_1,i_2,\dots ,i_p)$ such that
all multipliers in the left hand side of ~(\ref{dvaPTn}) are less
than $l$. This value equals
\begin{displaymath}\sum_{k=1}^{n}\sum_{\begin{array}{c} (m_1,\dots ,m_k)
\\ m\leq m_1,\dots ,m_k\leq l \\ m_1\cdot m_2 \dots m_k=l\cdot
(z_a +1)^{n-1}\end{array}}\binom{n}{k}
\overset{k}{\underset{j=1}{\prod}}\alpha_{m_j} .
\end{displaymath}
The number of sets $(p ; i_1 ,i_2,\dots ,i_n)$ such that one of
the multipliers in the left hand side of~(\ref{dvaPTn}) equals
$l$, is equal to $n\cdot\alpha_l$. Then $\alpha_l$ can be found
from the equality $A+n\cdot\alpha_l =C $.

Now let $m\leq z_a +1$. Then the cardinality of the equivalence
class of the relation $\sim_a$ is the least if in~(\ref{odynPTn})
all multipliers are equal to $m$, that is, this cardinality equals
$m^n$. Count the number of different equivalence classes
$\overline{x_0}$ of the relation $\sim_a$ of the cardinality
$m^n$. To make  all multipliers in~(\ref{odynPTn}) equal to $m$,
there should be $n_a (y_i)=m$ for all $i$. However $|\{t : n_a
(t)=m\}|=\alpha_m$. So every $y_i$ can be chosen in $\alpha_m$
ways. Since the meanings of $y_i$ for different $i$ are chosen
independently, there are $\alpha_{m}^{n}$ different $y=x_0 a$,
such that the corresponding class $\overline{x_0}$ contains $m^n$
elements. Hence by the relation $\sim_a$ we may define the number
$m$ and the meaning $\alpha_m$ of the first non zero component of
the type $(\alpha_1,\alpha_2,\dots ,\alpha_n)$ of $a$.

Components $\alpha_l$ for $l>m$ can be found recursively applying
the same arguments as above. Some changes include the following:
in ~(\ref{dvaPTn}) the right hand side must be substituted with
$l\cdot m^{n-1}$, and the number of sets $(p; i_1, i_2,\dots
,i_n)$, such that one of the multipliers in the left hand side
of~(\ref{dvaPTn}) equals $l$,
 and other $n-1$ multipliers equal $m$ is equal to $n\cdot\alpha_l \cdot\alpha_{m}^{n-1}$
 (if $z_a +1\neq l$  then $p=n$ ). Therefore $\alpha_l$
can be found from the equality
$A+n\cdot\alpha_l\cdot\alpha_{m}^{n-1}=C$.

Analogously we may find the type $(\beta_1,\beta_2,\dots
,\beta_n)$ of $b$ via the cardinalities of the equivalence classes
of the relation $\sim_b$. Since for isomorphic semigroups
$(\mathcal{PT}_n ,*_a)$ and $(\mathcal{PT}_n ,*_b)$ the number of
the equivalence classes of the same cardinality coincide, and the
values $\alpha_k$ and $\beta_k$, $k=1,\dots,n$ are defined
uniquely by the number of equivalence classes, for all $k$ we have
$\alpha_k =\beta_k$, that is, elements $a$ and $b$ have the same
types.

\textsl{Sufficiency.} Let elements $a$ and $b$ have the type
$(\alpha_1,\alpha_2,\dots, \alpha_n)$. Then there exist
permutations $\pi$ and $\tau$ in $\mathcal{S}_n $ such that
$b=\tau\alpha\pi$. The map $f : (\mathcal{PT}_n ,*_a)\rightarrow
(\mathcal{PT}_n ,*_b)$ such that $f(x)=\pi^{-1}x\tau^{-1}$ defines
isomorphism between $(\mathcal{PT}_n ,*_a)$ and $(\mathcal{PT}_n
,*_b)$. In fact $f$ is bijective and
\begin{align*} f(x *_a
y)=\pi^{-1}x *_a y\tau^{-1}=\pi^{-1}x\tau^{-1}\tau a\pi\pi^{-1}y
\tau^{-1}=\\
=\pi^{-1}x\tau^{-1}b\pi^{-1}y\tau^{-1}=f(x)*_b f(y).
\end{align*}
\end{proof}
\begin{corollary}
Let $p(k)$ denote the number of ways in which one can split
positive integer $k$ into non ordered sum of the natural integers.
Then there are $\overset{n}{\underset{k=0}{\sum}}p(k)$ pairwise
non-isomorphic semigroups obtained from $\mathcal{PT}_n$ by the
deformed multiplication.
\end{corollary}

\vspace{0.1cm} \noindent Department of Mechanics and
Mathematics,\\
Kiev Taras Shevchenko University,\\
 64, Volodymyrska st., 01033,
Kiev, Ukraine,\\ e-mail: {\em gtsyaputa\symbol{64}univ.kiev.ua}
$$ $$
\textit{Given to the editorial board} 30.09.2003
\end{document}